\documentclass[11pt]{amsart}

\newtheorem*{prop}{Proposition}
\newtheorem*{lem}{Lemma}
\newtheorem*{cor}{Corollary}

\theoremstyle{definition}
\newtheorem*{ack}{Acknowledgment}

\theoremstyle{remark}

\newcommand{\Irr}{\operatorname{Irr}}
\newcommand{\Id}{\operatorname{Id}}

\newcommand{\End}{\operatorname{End}}
\newcommand{\tr}{\operatorname{trace}}
\newcommand{\p}{\varphi}
\newcommand{\e}{\varepsilon}
\renewcommand{\d}{\delta}
\renewcommand{\l}{\lambda}
\renewcommand{\L}{\Lambda}
\renewcommand{\O}{\mathcal{O}}
\newcommand{\Z}{\mathbb{Z}}
\newcommand{\cl}{\operatorname{cl}}
\newcommand{\iso}{\overset{\sim}{\longrightarrow}}

\begin{document}


\title[Representations of Hopf Algebras]
	{On the Degrees of Irreducible Representations of Hopf Algebras}
\author{Martin Lorenz}
\address{Department of Mathematics, Temple University,
    Philadelphia, PA 19122-6094}
\email{lorenz@math.temple.edu}
\thanks{Research supported in part by NSF Grant DMS-9400643}

\begin{abstract}
Let $H$ denote a semisimple Hopf algebra over an algebraically closed field $k$ of
characteristic 0. We show that the degree of any irreducible representation of $H$
whose character belongs to the center of $H^*$ must be a divisor of $\dim_kH$.
\end{abstract}
\maketitle


\section{Introduction} \label{S:intro}

Let $H$ denote a semisimple Hopf algebra over an algebraically closed field $k$ of
characteristic 0. According to a famous conjecture of Kaplansky \cite{K}, the degrees 
of the irreducible representations of $H$ are expected to be divisors of $\dim_kH$.
This has been confirmed under the additional hypothesis that $H$ is
quasi-triangular \cite{EG}; see \cite{S1} for an elegant new proof.

In this note, we show by a very simple argument that the desired divisibility relation 
$\frac{\dim H}{\dim V}\in\Z$ also
holds for any irreducible $H$-module $V$ whose character $\chi_V$ belongs to
the center $Z(H^*)$ of $H^*$. Here, as usual, the \emph{character} of $V$ is the linear form
$\chi_V\in H^*$ that is defined by
$\chi_V(h) = \tr_{V/k}(h_V)$, where $h_V\in\End_k(V)$ is given by 
the action of $h\in H$. Thus: 

\begin{prop}
Let $H$ denote a semisimple Hopf algebra over an algebraically closed field $k$ of
characteristic 0. Then $\dim_kV$ divides $\dim_kH$ for any irreducible $H$-module $V$ 
such that $\chi_V\in Z(H^*)$.
\end{prop}

The work presented here was inspired by Y. Sommerh\"auser's lecture ``On central
character rings", delivered during the MSRI-workshop on Hopf algebras in October
1999.

The above notations will remain in effect throughout. In addition, $\Irr(H)$ will
denote a fixed representative set of the isomorphism classes of irreducible (left) 
$H$-modules and $C(H)$ is the \emph{character algebra}, that is, the $k$-span in 
$H^*$ of all characters $\chi_V$ of $H$-modules $V$. 
In general, our notation for Hopf algebras
follows \cite{Mont}.


\section{Some Preliminaries} \label{S:prelim}

Our hypotheses entail that the dual Hopf algebra $H^*$ is semisimple as well; 
see \cite{Mont}.
Recall that $H$ acts on $H^*$ via $\langle h\p,h'\rangle =\langle \p,h'h\rangle$
for $h,h'\in H$ and $\p\in H^*$. Dually, $H^*$ acts on $H$. Some particulars
of these actions are explained in the following standard lemma originally due to Masuoka
\cite{M}; see also \cite[p.~44]{S2}. 

We let $\l\in H^*$ denote the integral 
with $\langle\l,1\rangle=1$ and $\L\in H$ the integral with
$\langle\l,\L\rangle=1$; cf. \cite[p.~26]{S2}. 
Moreover, $S^*$ denotes the antipode of $H^*$.

\begin{lem}\label{lem1}
Let $V\in\Irr(H)$ and let $e_V\in Z(H)$ denote the corresponding centrally 
primitive idempotent of $H$, acting as $\Id_V$ on $V$ and as $0_W$ on all
$W\in\Irr(H)\setminus\{V\}$. Then
$$
\tfrac{\dim H}{\dim V}e_V = (S^*\chi_V)\L
\qquad\text{and}\qquad
\tfrac{\dim H}{\dim V}e_V\l = \chi_V\ .
$$
\end{lem}

\begin{proof}
The first formula 
is Corollary 4.6 in \cite{S2} and the second one is
established in the proof of \cite[Lemma (b)]{M}; cf. also 
\cite[proof of Proposition 4.5]{S2}. 
\end{proof}

\begin{cor}
For every idempotent $\d=\d^2\in Z(H^*)$, the element $\d\L\in C(H^*)\subseteq H$
is a character of some $H^*$-module.
\end{cor}

\begin{proof}
We will apply the Lemma with the roles of $H$ and $H^*$ interchanged. Since
$\langle\e,\L\rangle=\dim H$ holds for the counit $\e=1_{H^*}$ (e.g., \cite[Prop.~3.4]{S2}),
this requires replacing $\L$ by $\L'=\frac{1}{\dim H}\L$.

First assume that $\d=\d_M$ is the centrally primitive idempotent of $H^*$ corresponding
to some irreducible $H^*$-module $M$. Then, by 
the Lemma, 
$\frac{\dim H^*}{\dim M}\d_M\L'=\chi_M\in C(H^*)$. Thus,
$$
\d_M\L=(\dim M)\chi_M=\chi_{M^{\dim M}} \ .
$$
In general, $\d=\sum_M\d_M$ for certain irreducible $H^*$-modules $M$,
and so $\d\L=\chi_N$ for the $H^*$-module $N=\oplus_M M^{\dim M}$. 
\end{proof}


\section{The Proof} \label{S:proofprop}

We are now ready to give the proof of the Proposition. So let $V\in\Irr(H)$
be given with $\chi_V\in Z(H^*)$. Our goal is to show that
$\frac{\dim H}{\dim V}\in \O$, where $\O$ denotes the
ring of algebraic integers in $k$. Since 
$\frac{\dim H}{\dim V}e_V = (S^*\chi_V)\L$, by the Lemma, it suffices to show
that $(S^*\chi_V)\L$ is integral over $\Z$.

To this end, let $\d_1,\ldots,\d_r$ be the distinct primitive idempotents of
$Z(H^*)$ and let $f_i:Z(H^*)\to k$ denote the corresponding characters; so
$\zeta=\sum_i f_i(\zeta)\d_i$ holds for all $\zeta\in Z(H^*)$. In particular,
since $S^*\chi_V\in Z(H^*)$, we have
$$
(S^*\chi_V)\L=\sum_i f_i(S^*\chi_V)\d_i\L \ .
$$
By the Corollary, each $\d_i\L$ is the character of some (irreducible)
$H^*$-module; so $\d_i\L$ belongs to the \emph{Grothendieck ring}
$G_0(H^*)=\bigoplus_{M\in\Irr(H^*)}\Z\chi_M\subseteq C(H^*)$. Every element of
$G_0(H^*)$ satisfies a monic polynomial over $\Z$. Analogously, the
character $S^*\chi_V=\chi_{V^*}\in G_0(H)$ satisfies a monic polynomial over $\Z$, and
so $f_i(S^*\chi_V)\in\O$ holds for all $i$. Consequently,
$$
(S^*\chi_V)\L \in G_0(H^*)\otimes\O \subseteq C(H^*)\ ,
$$
which entails that $(S^*\chi_V)\L$ is integral over $\Z$, as desired.


\section{Concluding Remarks} \label{S:remarks}

Let $(\,.\,)^{\cl}$ denote the set of elements in the ring 
in question that are integral over $\Z$; e.g., 
$$
Z(H)^{\cl}=\bigoplus_{V\in\Irr(H)}\O e_V\ .
$$
Consider the isomorphism $f\colon H^*\iso H$, $f(\p)=\p\L$.
By the Lemma, $f$ restricts to isomorphisms $C(H)\iso Z(H)$ and
$Z(H^*)\iso C(H^*)$. The essence of the above proof is that, in fact,
$$
Z(H^*)^{\cl} \stackrel{f}{\longrightarrow} 
G_0(H^*)\otimes\O \subseteq C(H^*)^{\cl}\ ,
$$
while Kaplansky's conjecture is equivalent with 
$$
G_0(H) \stackrel{f}{\longrightarrow} Z(H)^{\cl}\ .
$$
It is tempting to try and consolidate the Etingof-Gelaki result 
for quasi-triangular Hopf algebras
and our Proposition by at least showing that $f$ maps the \emph{center}
of $G_0(H)$ to $Z(H)^{\cl}$. 


\begin{ack} The work on this note was started while I was
attending the Noncommutative Algebra program at MSRI in the Fall of 1999. 
I would like to thank the organizers of this program, in particular S. Montgomery, 
and MSRI staff for their hospitality and support.
I would also like to thank Temple University for granting me
a research leave during the Fall Semester 1999.
\end{ack}


\end{document}